\documentclass[3p,,preprint,12pt]{elsarticle}

\usepackage{amssymb}
\usepackage{endnotes}
\usepackage{enumerate}
\usepackage{footmisc}
\usepackage{amssymb}
\usepackage{amsmath}
\usepackage{soul}
\usepackage{bm}
\usepackage{float}

\usepackage{multirow}
\usepackage{makerobust}
\MakeRobustCommand\overrightarrow
\MakeRobustCommand\overleftarrow
\usepackage{pdflscape}
\usepackage{endnotes}
\usepackage{amsmath}
\usepackage{graphicx}
\usepackage{caption}
\usepackage{subcaption}
\usepackage{endnotes}
\usepackage{pgf}
\usepackage{epstopdf}
\usepackage{float}
\usepackage{appendix}
\usepackage{chngcntr}
\usepackage{apptools}

\newtheorem{theorem}{Theorem}
\newtheorem{lemma}{Lemma}
\newtheorem{proposition}{Proposition}
\newtheorem{corollary}{Corollary}

\newtheorem{definition}{Definition}

\newtheorem{property}{Property}

\begin{document}

\begin{frontmatter}

\title{A Marginal Analysis Framework to Incorporate the Externality Effect of Ordering Perishables}

\author[afe5e335e84b9]{Katsunobu Sasanuma\corref{c-69426cb9627c}}

\ead{katsunobu.sasanuma@stonybrook.edu}\cortext[c-69426cb9627c]{Corresponding author.}

\author[afe5e335e84b9]{Mohammad Delasay}

\ead{mohammad.delasay@stonybrook.edu}

\author[afe5e335e84b9]{Christine Pitocco}

\ead{christine.pitocco@stonybrook.edu}

\author[a67c87fff5111]{Alan Scheller-Wolf}

\ead{awolf@andrew.cmu.edu}

\author[afe5e335e84b9]{Thomas Sexton}

\ead{thomas.sexton@stonybrook.edu }

\address[afe5e335e84b9]{College of Business\unskip, Stony Brook University\unskip, Stony Brook\unskip, NY 11794\unskip, USA}

\address[a67c87fff5111]{Tepper School of Business\unskip, Carnegie Mellon University\unskip, Pittsburgh\unskip, PA 15213\unskip, USA}

\begin{abstract}
Finding the optimal policy for multi-period perishable inventory systems requires solving computationally-expensive stochastic dynamic programs (DP). To avoid the difficulty of solving DP models, we propose a framework that uses an externality term to capture the long-term impact of ordering decisions on the average cost over an infinite horizon. By approximating the externality term, we yield a tractable approximate optimality condition, which is solved through standard marginal analysis. The resulted policy is near-optimal in long-run average cost and ordering decisions.
\end{abstract}

\begin{keyword}
Perishable inventory \sep marginal analysis \sep externality \sep constant base-stock policy
\end{keyword}

\end{frontmatter}

\section{Introduction}
The exact analysis of perishable inventory models using stochastic dynamic programming (DP) is computationally expensive, rendering it intractable for models with large state spaces. As noted in~\citet{karaesmen2011managing}, ``The policy structures outlined in~\citet{fries1975optimal} and~\citet{nahmias1975optimal} are quite complex; perishability destroys the simple base-stock structure of optimal policies for discrete review models without fixed order costs in the absence of perishability.''
Furthermore, the DP approach does not provide any insight into the form of the inventory-dependent (which we refer to as \emph{state-dependent}) optimal policy. Many researchers have thus sought effective heuristic methods \citep[for comprehensive reviews, see, e.g.,][]{karaesmen2011managing,baron2011managing}. Among these heuristics, the constant base-stock (CBS) policy---despite its simplicity---has been shown to be an excellent alternative to the optimal state-dependent policies~\citep{huh2009asymptotic,bijvank2014robustness, bu2020asymptotic}. Many state-dependent policies have also been proposed, among which two approaches have received increased attention, namely~\emph{$L^\natural$-convexity}~\citep[e.g.,][]{zipkin2008structure,huh2010optimal,chen2014coordinating} and the \emph{marginal cost accounting scheme}~\citep[e.g.,][]{levi20082, truong2014approximation, chao2015approximation, zhang2016approximation}. The marginal cost accounting scheme utilizes marginal analysis, which provides an efficient algorithm for perishable inventory models. The key of the marginal cost accounting scheme is to develop an effective cost-balancing technique for the specific model under consideration, which is often not straightforward to identify.

We develop a marginal analysis framework that incorporates the \emph{externality} effect---the indirect long-term impact of ordering decisions---on the average cost of a perishable inventory system. To our knowledge, the inclusion of the externality effect in a marginal analysis framework has not been employed in the inventory management literature, though it has been widely implemented to study economic concepts, including congestion pricing~\citep[e.g.,][]{vickrey1969congestion,larson2010}. Using this framework, we derive an approximate optimality condition for the general state-dependent policy. This optimality condition is a recursive equation, which unfortunately is difficult to solve. However, by utilizing the properties of the CBS policy, we can reduce the externality effect into a fixed cost or benefit, representing the marginal external cost. Thus, we convert the original complex exact optimality condition into a simple approximate optimality condition in which only a single order amount for a given initial inventory level is involved. This single-decision condition is \emph{almost} identical to the optimality condition for the newsvendor model, and hence, is easy to solve, for example, in a spreadsheet. This approach provides near-optimal solutions both with respect to the average cost and the individual order amounts. In addition, our approach provides insight into the state-dependent characteristics of near-optimal ordering policies.

There is an abundance of near-optimal heuristics for perishable inventory systems in the literature. In this paper, our primary contribution is not to add one more element to this list, but rather to provide a general framework to convert numerically intractable multi-decision stochastic dynamic inventory models to tractable single-decision models. Our framework is motivated by density functional theory~\cite{hohenberg1964inhomogeneous} and its local density approximation~\cite{kohn1965self}---the most popular and successful methods in computational physics and chemistry to convert multi-body problems to single-body problems. We showcase the application of our framework on one of the classic perishable inventory models, but we believe it has the potential to be applied to other inventory models as well.

\vspace{-0.3cm}
\section{General Formulation}
\label{sec:general for}
In this section we describe our marginal analysis framework we use to derive the optimality condition of the state-dependent policy for a general infinite-horizon inventory system with a single perishable product. (In~\S\ref{sec:application}, we illustrate how to apply it to the model introduced in~\cite{nahmias1975optimal}.) In an infinite-horizon single-product perishable inventory system, the initial inventory is reviewed in each period, a new order is placed, demand is fulfilled, and  perished products are discarded. Let~$m(\ge 1)$ be the product lifetime and~$x_i$ be the number of units with~$i \in [1,m]$ remaining periods of lifetime. Then, the \emph{initial inventory} and \emph{order amount} are represented by~$x_1, \dots, x_{m-1}$ and~$x_m$, respectively. For ease of representation, let~$x^i=\sum_{j=1}^{i} x_j$ be the inventory level with remaining lifetime of at most~$i$ periods and~$x=x^{m-1}$ be the total initial inventory. Similarly, let~$\mathbf{x}^i=(x_1,...,x_i)$ be the inventory vector with the remaining lifetime of at most~$i$ periods and~$\mathbf{x}=\mathbf{x}^{m-1}$ be the total initial inventory vector. For notational convenience, let~$x^0=\mathbf{x^0}=0$. Let~$\Omega=\mathbb{R}_{\geq0}^{m-1}$, a set of non-negative real numbers in an~$(m-1)$-dimensional vector space; then~$\mathbf{x} \in \Omega, \forall \mathbf{x}$.

A stationary ordering policy may be characterized by its order-up-to level~$q(\mathbf{x})$, which is a scalar valued function of the initial inventory vector~$\mathbf{x}$. With a slight abuse of notation, we denote~$q$ to represent either the order-up-to level~$q(\mathbf{x})$ for a particular~$\mathbf{x}$ or the policy~$q(\cdot)$, a function of~$\mathbf{x}$, distinguishing between the two when necessary. When implementing the policy~$q$, the order amount at the beginning of each period becomes~$x_m=\max\{q(\mathbf{x})-x,0\}$, or simply~$x_m=[q-x]^+$.

We propose a stationary model of this infinite-horizon problem based on the ensemble-average cost (taken over the initial inventory distribution), instead of its stochastic dynamic program (DP) model, which is known to be computationally difficult to solve due to the curse of dimensionality and dependence of decisions among different periods. In our stationary model, the complexity of tracking inventory levels in infinite time periods is incorporated into the initial inventory distribution. The infinite-horizon DP model and the stationary model represent the same average total cost; DP calculates the time-average cost, and the stationary model calculates the ensemble-average cost.

When demand is independent and identically distributed (i.i.d.), we can define each period's initial inventory~$\mathbf{X}$ as a non-negative \emph{random vector} following a stationary distribution~$f_{\mathbf{X}}^q(\cdot)$ given the policy~$q$. Let~$L(q, \mathbf{x})$ be the one-period cost associated with the single ordering decision~$q$ when the initial inventory~$\mathbf{x}$ is observed at the beginning of the period. Then, the average total cost of the stationary model follows:
\begin{small}
\begin{align}
L(q) = \mathbb{E}_{\mathbf{X}}[L(q,\mathbf{X})]= \int_\Omega L(q(\mathbf{k}),\mathbf{k}) f_{\mathbf{X}}^q(\mathbf{k})d\mathbf{k},
\label{eq:expected total cost}
\end{align}
\end{small}
\vspace{-0.4cm}

\noindent which is a functional of the policy~$q=q(\mathbf{x}), \forall \mathbf{x} \in \Omega$. Let~$q^*$ be the minimizer of~$L(q)$ (i.e,~$q^*$ is the optimal inventory policy). When~$L(q)$ is convex (which is the case for many inventory models including the perishable inventory model in~\S\ref{sec:application}),~$q^{\ast}$ satisfies the following optimal functional derivative condition: 
\begin{small}
\begin{align}
\dfrac{\delta L(q)}{\delta q(\mathbf{x})}&=\dfrac{\partial L(q(\mathbf{x}),\mathbf{x})}{\partial q(\mathbf{x})} f_{\mathbf{X}}^{q}(\mathbf{x})+ \int_\Omega L(q(\mathbf{k}),\mathbf{k})\dfrac{\partial f_{\mathbf{X}}^q(\mathbf{k})}{\partial q(\mathbf{x})}d\mathbf{k} = 0, \quad \forall \mathbf{x} \in \Omega.
\label{eq:stationary optimality condition}
\end{align}
\end{small}
\vspace{-0.2cm}

The derivation of~\eqref{eq:stationary optimality condition} is motivated by the Kohn-Sham approach to reduce the dimensionality of multi-body problems in Physics~\cite{kohn1965self}; the details of the derivations and proofs appear in the online appendix. The optimality condition~\eqref{eq:stationary optimality condition} for the stationary model has two components. When the order-up-to level (i.e., the policy) changes from~$q(\mathbf{x})$ to~$q(\mathbf{x})+ \delta q(\mathbf{x})$ for the initial inventory~$\mathbf{x}$: (i) The first term is the contribution of this policy change to the average total cost~$L(q)$, assuming the initial inventory distribution remains the same; and (ii) The second term is the contribution of the policy change to the average total cost~$L(q)$ due to the \emph{change in the initial inventory distribution}. This second term, which we refer to as \emph{externality}, captures \emph{the long-term impact of ordering decisions}, since an equilibrium inventory distribution is reached only after infinitely many periods.

The externality term is the main source of complexity in the exact optimality condition~\eqref{eq:stationary optimality condition}; specifically, it is difficult to evaluate the function~${\partial f_{\mathbf{X}}^q(\mathbf{k})}/{\partial q(\mathbf{x})}$ representing the impact of the policy change at~$\mathbf{x}$ on the distribution at all~$\mathbf{k}\in \Omega$. To resolve this complexity, we can approximate the externality term using any simple and reasonably good policy~$\widetilde{q}$\:: Specifically, we replace~${\partial f_{\mathbf{X}}^{q}(\mathbf{k})}/{\partial {q}(\mathbf{x})}$ in the externality term with~${\partial f_{\mathbf{X}}^{\widetilde{q}}(\mathbf{k})}/{\partial {\widetilde{q}}(\mathbf{x})}$ along with some necessary modifications due to normalization. Adopting the idea of the local density approximation~\cite{kohn1965self}, we utilize the CBS policy, which is simple to optimize and known to be a reasonably good policy for many inventory models. By using the optimal CBS policy~$q_c^{\ast}$ instead of the optimal policy~$q^{\ast}$, we simplify the externality term in~\eqref{eq:stationary optimality condition} to~\eqref{eq: approximate externality}, in which~$V_{ex}(q_c)$ follows~\eqref{eq: externality term}. We expect this to be a good approximation because: (i) The expected change in the one-period cost originating from the change of the initial inventory distribution is conserved:~$L(q^{\ast}(\mathbf{x}),\mathbf{x}) \delta f^{q^{\ast}}_{\mathbf{X}}(\mathbf{x})\approx L(q_c^{\ast},\mathbf{x}) \delta f^{q^{\ast}_c}_{\mathbf{X}}(\mathbf{x})$; (ii) The expected change in the order amount is conserved:~$f^{q^{\ast}}_{\mathbf{X}}(\mathbf{x})\delta q^{\ast}(\mathbf{x}) \approx \delta q^{\ast}_c$ (or equivalently,~$\partial q_c /\partial q(\mathbf{x})\big|_{q=q^{\ast}} \approx f^{q^{\ast}}_{\mathbf{X}}(\mathbf{x})$). Therefore 
\begin{small}
\begin{gather}
\int_{\Omega} L(q^{\ast}(\mathbf{k}),\mathbf{k}) \left.{\dfrac{\partial f^q_{\mathbf{X}}(\mathbf{k})}{\partial q(\mathbf{x})}} \right|_{q=q^{\ast}} d\mathbf{k} \approx \int_{\Omega} L(q_c^{\ast}(\mathbf{k}),\mathbf{k}) \left.{\dfrac{\partial f^{q_c}_{\mathbf{X}}(\mathbf{k})}{\partial q_c}}\right|_{q_c=q_c^{\ast}} \left.{\dfrac{\partial q_c}{\partial q(\mathbf{x})}}\right|_{q=q^{\ast}}d\mathbf{k} \approx V_{ex}(q^{\ast}_c) f^{q^{\ast}}_{\mathbf{X}}(\mathbf{x}),
\label{eq: approximate externality}\\
 \qquad \qquad \qquad V_{ex}(q_c)=\int_{\Omega} L(q_c,\mathbf{k}) \dfrac{\partial f^{q_c}_{\mathbf{X}}(\mathbf{k})}{\partial q_c}d\mathbf{k}.
\label{eq: externality term}
\end{gather}
\end{small}

\vspace{-0.4cm} Combining~\eqref{eq:stationary optimality condition} and~\eqref{eq: approximate externality}, we derive the approximate optimality condition~\eqref{eq:optimality condition} for the stationary model, conditioned on~$\mathbf{x}$ being recurrent (i.e.,~$f^q_{\mathbf{X}}(\mathbf{x})>0$):
\begin{small}
\begin{align}
    \left( \dfrac{\partial L(q,\mathbf{x})}{\partial q(\mathbf{x})} + V_{ex}(q^{\ast}_c) \right) f^q_{\mathbf{X}}(\mathbf{x})=0, \forall \mathbf{x} \in \Omega \implies \dfrac{\partial L(q,\mathbf{x})}{\partial q(\mathbf{x})}+ V_{ex}(q^{\ast}_c)=0.\label{eq:optimality condition}
\end{align}
\end{small}

\vspace{-0.2cm} Similar to~\eqref{eq:stationary optimality condition}, the optimality condition~\eqref{eq:optimality condition} has two components. We refer to the first term~${\partial L(q,\mathbf{x})}/{\partial q(\mathbf{x})}$ as the \emph{marginal internal cost} ($M\!IC$) and the second term~$V_{ex}(q^{\ast}_c)$ as the \emph{marginal external cost} ($M\!EC$), which is a constant since it is independent of~$\mathbf{x}$ under the CBS approximation. Without the $M\!EC$ term,~\eqref{eq:optimality condition} reduces to the optimality condition of a standard single-decision inventory model, which is easy to solve. But the $M\!EC$ term does not increase the computational complexity of solving~\eqref{eq:optimality condition} as it is simply a constant. Nevertheless, it plays an important role in minimizing the average cost. By Solving~\eqref{eq:optimality condition}, we obtain the approximate optimal policy~$q_h^{\ast}$, which is a state-dependent policy (due to the~$M\!IC$ term) like the optimal policy.
\vspace{-0.3cm}
\section{Applying the Framework for a Perishable Inventory Model}
\label{sec:application}
In this section we showcase how our framework, described in~\S\ref{sec:general for}, can be used to analyze the classic model by~\citet{nahmias1975optimal}: A periodic-review perishable inventory model with lost sales, fixed product lifetime, no lead time, i.i.d.~demand (denoted by~$D$), and a FIFO issuance policy. The four cost parameters include: purchase ($c$ per unit), holding ($h \geq 0$ per period per unit), shortage ($r$ per unit), and wastage ($\theta$ per unit). The optimal policy for this model is obtained using DP in~\cite{nahmias1975optimal}. In this section, we show how to derive the approximate optimality condition~\eqref{eq:optimality condition} for this model, and then, we evaluate its accuracy.

To derive the one-period cost~$L(q,\mathbf{x})$, we evaluate the costs associated with a single ordering decision at the beginning of period~1, considering that the holding and shortage costs are incurred in period~1 and wastage cost is incurred in period~$m$. Let the random variable~$D^i(\mathbf{x}^{i-1})$ represent the~$i$-period \emph{effective demand}, i.e., the total outflow through demand and wastage from periods 1 to~$i$ (excluding the wastage in period~$i$). Denoting the random variables for demand and wastage in period~$i$ by~$D_i$ and~$R_i$, respectively, we can define~$D^i(\mathbf{x}^{i-1})$ recursively:
\vspace{-0.3cm}
\begin{small}
\begin{align}
D^{i}(\mathbf{x}^{i-1}) =
  \begin{cases}
   D^{i-1}(\mathbf{x}^{i-2})+R_{i-1}+D_{i} & \text{if }  i = 2, \dots, m,\\
   D_1  & \text{if } i=1,
  \end{cases}  \label{eq: Di}\\
R_{i-1} =
  \begin{cases}
   [x^{i-1} - D^{i-1}(\mathbf{x}^{i-2})]^+ & \text{if } i = 2, \dots, m,\\
   0  & \text{if } i=1.
  \end{cases}
  \label{eq: Ri}
\end{align}
\end{small}
\vspace{-0.4cm}

Note that~$D^{m}(\mathbf{x}^{m-1})=D^m(\mathbf{x})$. We denote the number of units being held, in shortage, and wasted under policy~$q$ as~$[q - D]^+$,~$[D - q]^+$, and~$[q - D^m(\mathbf{x})]^+$, and we express their corresponding expectations by~$n_h(q)$,~$n_s(q)$, and~$n_w(q)$, respectively. Then, we can represent the expected one-period cost~$L(q,\mathbf{x})$ and the average total cost~$L(q)$:\footnote{For ease of exposition, we incorporate the purchase cost when the unit is either sold or perished; i.e.,~$r-c>0$ and~$\theta+c>0$ represent the shortage (understocking) and the wastage (overstocking) costs.}
\begin{small}
\begin{gather}
L(q,\mathbf{x})=h \mathbb{E}_D[q - D]^+ + (r-c) \mathbb{E}_D[D - q]^+ + (\theta+c) \mathbb{E}_{D^m(\mathbf{x})}[q - D^m(\mathbf{x})]^+.\\
L(q) =h  \underbrace{\mathbb{E}_{\mathbf{X}} \mathbb{E}_D[q(\mathbf{X}) - D]^+}_{n_h(q)} + (r-c) \underbrace{\mathbb{E}_{\mathbf{X}} \mathbb{E}_D[D - q(\mathbf{X})]^+}_{n_s(q)} + (\theta+c) \underbrace{\mathbb{E}_{D^m(\mathbf{X})}[q(\mathbf{X}) - D^m(\mathbf{X})]^{+}}_{n_w(q)}.
\label{eq:expected total cost2}
\end{gather}
\end{small}
\vspace{-0.5cm}

\paragraph{\textbf{Marginal Internal Cost ($M\!IC$)}} To evaluate the $M\!IC$ term, we need~$L(q,\mathbf{x})$, which in turn depends on the~$m$-period effective demand~$D^m(\mathbf{x})$. Proposition~\ref{prop:convolution formula} shows how to obtain its cumulative distribution function (CDF).
\begin{proposition}
\label{prop:convolution formula}
The CDF of~$D^m(\mathbf{x})$ is obtained by applying the following recursively:
\begin{small}
\begin{align}
\label{eq:convolution formula}
F_{D^{i+1} (\mathbf{x}^{i})}(z)=
	\begin{cases}
	\int_{\xi=0}^{z-x^i}  F_{D^{i} (\mathbf{x}^{i-1})} (z-\xi) f_{D_{i+1}}(\xi) d \xi & \quad \text{if } z>x^i, i=1, \dots, m-1,\\
	0 & \quad \text{if } z \le x^i, i=1, \dots, m-1,\\
	F_{D}(z) & \quad \text{if } i=0.
	\end{cases}
\end{align}
\end{small}
\end{proposition}

We can represent~$L(q,\mathbf{x})$ and its partial derivative (i.e., $M\!IC$) as~\eqref{eq:total cost given initial inventory} and~\eqref{eq:marginal internal cost}, respectively. We can confirm that~$L(q,\mathbf{x})$ is strictly convex as~${\partial^2 L(q,\mathbf{x})}/{\partial q(\mathbf{x})^2}>0,\forall q \in [0,\infty)$.
\vspace{-0.5cm}
\begin{small}
\begin{gather}
L(q,\mathbf{x})=
h\int_0^{q} (q-z)f_D(z)dz + (r-c) \int_{q}^\infty (z-q)f_D(z)dz + (\theta+c) \int_{0}^{q}(q-z)f_{D^m(\mathbf{x})}(z)dz.
\label{eq:total cost given initial inventory}\\
M\!IC: \dfrac{\partial L(q,\mathbf{x})}{\partial q(\mathbf{x})}=-(h+r-c)\bar{F}_D(q)+(\theta+c)F_{D^m(\mathbf{x})}(q)+h,
\label{eq:marginal internal cost}
\end{gather}
\end{small}
\noindent where~$\bar{F}_D(q)$ is the complementary CDF of the demand distribution.

\paragraph{\textbf{Marginal External Cost ($M\!EC$)}} To evaluate the $M\!EC$ term~$V_{ex}(q_c)$, we substitute the expression for~$L(q_c,\mathbf{x})$ from~\eqref{eq:total cost given initial inventory} into~\eqref{eq: externality term}. The first and second terms in~\eqref{eq:total cost given initial inventory} do not contribute to the $M\!EC$ term as they do not depend on~$\mathbf{x}$.\footnote{This is because~$\int_{\Omega} \dfrac{\partial f^{q_c}_{\mathbf{X}}(\mathbf{x})}{\partial q_c}d\mathbf{x}= \dfrac{\partial [\int_{\Omega} f^{q_c}_{\mathbf{X}}(\mathbf{x})d\mathbf{x}]}{\partial q_c}= \dfrac{\partial 1}{\partial q_c}= 0$, where we use Leibniz's rule.}  As a result, we derive~\eqref{eq:marginal external cost2}. Let~$\mathbf{X}^{q_c}$ denote the initial inventory random vector under the CBS policy~$q_c$. We can compute~$w_{ex}(q_c)$ following~\eqref{eq:marginal external cost} by discretizing~$q_c$ with step size~$\Delta$ and evaluating the difference between two expectations. According to Proposition~\ref{prop:external cost},~$w_{ex}(q_c)$ is bounded.
\begin{small}
\begin{gather}
M\!EC: V_{ex}(q_c)=(\theta+c)\int_\Omega \int_{0}^{q_c}(q_c -z)f_{D^m(\mathbf{k})}(z)dz \dfrac{\partial f_{\mathbf{X}}^{q_c}(\mathbf{k})}{\partial q_c} d \mathbf{k}=(\theta+c) w_{ex}(q_c)\label{eq:marginal external cost2},\\
w_{ex}(q_c)=\frac{1}{\Delta}\Big( \underbrace{\mathbb{E}[q_c - D^m(\mathbf{X}^{q_c+\Delta})]^{+}}_{n_w^\Delta(q_c)}-\underbrace{\mathbb{E}[q_c - D^m(\mathbf{X}^{q_c})]^{+}}_{n_w(q_c)}\Big).
\label{eq:marginal external cost}
\end{gather}
\end{small}
\vspace{-.6cm}

\begin{proposition}
\label{prop:external cost}
The externality is negative and bounded; i.e.~$-1 < w_{ex}(q_c) \leq 0, \forall q_c \geq 0.$
\end{proposition}

By replacing the $M\!IC$ term~\eqref{eq:marginal internal cost} and the $M\!EC$ term~\eqref{eq:marginal external cost2} in the approximate optimality condition~\eqref{eq:optimality condition}, we obtain the approximate optimality condition as follows:
\begin{small}
\begin{align}
\label{eq:optimality3}
(\theta+c) F_{D^m(\mathbf{x})}(q)=(h+r-c)\bar{F}_D(q)-h-(\theta+c) w_{ex}(q^{\ast}_c).
\end{align}
\end{small}
\vspace{-0.5cm}

The solution to~\eqref{eq:optimality3} is unique (Proposition~\ref{thm:existence}). We can thus find the approximate optimal policy~$q_h^{\ast}(\mathbf{x}), \forall \mathbf{x} \in \Omega$ numerically. Based on~\eqref{eq:optimality3}, the approximate optimal policy approaches CBS as~$h$ or~$r$ grow large; the same patterns hold for the actual optimal policy~\citep{huh2009asymptotic}. Also, the optimal policy is asymptotically CBS when demand variability decreases or~$m$ increases~\citep{bu2020asymptotic}; the same patterns hold for our approximate optimal policy, following Proposition~\ref{thm:CBS condition}.
\begin{proposition}
\label{thm:existence}
There exists a unique finite order-up-to level (approximate optimal policy)~$q_h^\ast(\mathbf{x})$ satisfying the optimality condition~\eqref{eq:optimality3} for any initial inventory vector~$\mathbf{x} \in \Omega$.
\end{proposition}
\begin{proposition}
\label{thm:CBS condition}
The approximate policy approaches CBS if and only if~$F_{D^m(\mathbf{0})}(q_h^{\ast}(\mathbf{0})) \to 0$.
\end{proposition}

We compare the average total costs under the approximate optimal policy (using the algorithm laid out in Table~\ref{tab:algorithm}) and the optimal policy (using the DP algorithm described in~\cite{nahmias1975optimal}), for~$c=0$ and other parameters as specified in Table~\ref{table:experiments1} under exponential and Poisson demands with mean~10. Note that the optimal policy is the same for any combination of the parameters that result in the same~$r-c$ and~$\theta+c$. As~$m$ grows, the optimal policy approaches CBS; therefore, the relative gap~$G$ between the optimal and our approximate policies is going to be more stark (if such a gap exists) for small~$m$ values, where the optimal policy is highly state-dependent. According to our numerical experiments, the average and maximum cost deviations between the approximate and optimal policies are very small---around~0.05\% and~0.34\%, respectively.
\begin{table}[h]
    \centering
\footnotesize
        \caption{Algorithm}
    \label{tab:algorithm}
    \begin{tabular}{l}
         \hline
          Pre-processing (performed for each combination of~$m$ and demand distribution):\\
         $\qquad$ Derive~$F_{D^m(\mathbf{x})}(q)$ (eq.~\eqref{eq:convolution formula}).\\
         $\qquad$ Discretize~$q_c$ and~$\mathbf{x}$ for continuous distributions. For each~$q_c$,  simulate a system with CBS policy\\ $\qquad$ and evaluate~$(n_{h}(q_c),n_{s}(q_c),n_{w}(q_c),n^\Delta_{w}(q_c))$ (eqs.~\eqref{eq:expected total cost2} and~\eqref{eq:marginal external cost}).\\
         Marginal analysis (performed for each combination of~$c$,~$h$,~$r$, and~$\theta$):\\
  	  $\qquad$ For each value of~$q_c$, evaluate~$L(q_c)$ (eq.~\eqref{eq:expected total cost2}).\\
  	  $\qquad$ Find~$q_c^{\ast}=arg\,min_{q_c} L(q_c)$ and evaluate $w_{ex}(q_c^{\ast})$ (eq.~\eqref{eq:marginal external cost}).\\
    $\qquad$ For each~$\mathbf{x} \in \Omega$, conduct marginal analysis to determine~$q_h^{\ast}(\mathbf{x})$ (eq.~\eqref{eq:optimality3}).\\
    $\qquad$ Find the order amount~$x_m=[q^{\ast}_h(\mathbf{x})-x]^+$.\\
    \hline
    \end{tabular}
\end{table}
 \begin{table}[h]
\caption{The comparison between the DP and approximate policies}
\centering
\footnotesize
\begin{tabular}{  l | c  c c| c  c  c| c c c| c c c }
\hline
	 & \multicolumn{6}{c|}{Exponential Demand} & \multicolumn{6}{c}{Poisson Demand}  \\
	 & \multicolumn{3}{c|}{$m=2$}      & \multicolumn{3}{c|}{$m=3$}  & \multicolumn{3}{c|}{$m=2$}  & \multicolumn{3}{c}{$m=3$}   \\
	$h, r, \theta$ & DP & $G\%$ & MAD & DP & $G\%$ & MAD & DP & $G\%$ & MAD & DP & $G\%$ &  MAD \\\hline
	0, 5, 5 & 19.84 & 0.04 & 0.28 & 12.14 & 0.34 & 0.62 & 1.47 & 0.14 & 0.13 & 0.13 & 0.27 & 0.00 \\
	0, 5, 10 & 25.40 & 0.06 & 0.26 & 16.05 & 0.07 & 0.26 & 2.09 & 0.11 & 0.20 & 0.19 & 0.04 & 0.11  \\
	0, 5, 20 & 30.74 & 0.02 & 0.14 & 20.24 & 0.02 & 0.33 & 2.92 & 0.08 & 0.07 & 0.26 & 0.00 & 0.00  \\
	0, 8, 7 & 30.06 & 0.05 & 0.29 & 18.31 & 0.20 & 0.48 & 2.16 & 0.00 & 0.06 & 0.19 & 0.00 & 0.00  \\
	0, 10, 5 & 29.19 & 0.09 & 0.38 & 17.49 & 0.17 & 0.46 & 1.95 & 0.00 & 0.06 & 0.17 & 0.00 & 0.00 \\
	1, 5, 5 & 25.39 & 0.01 & 0.14 & 20.88 & 0.00 & 0.05 & \underline{5.26} & 0.27 & 0.43 & \underline{4.93} & 0.00 & 0.00  \\
	1, 5, 10 & 28.93 & 0.01 &0.14 & 22.69 & 0.00 &0.04 & 5.52 & 0.00 &0.00 & \underline{4.93} & 0.00 & 0.00  \\
	1, 5, 20 & 32.81 & 0.00 &0.05 & 25.03 & 0.01 &0.10 & 5.88 & 0.00 & 0.00& \underline{4.94} & 0.00 &  0.00 \\
	1, 8, 7 & 36.51 & 0.02 &0.17 & 28.38 & 0.02 & 0.16& 6.36 & 0.00 &0.00 & \underline{5.68} & 0.00 &  0.00 \\
	1, 10, 5 & 38.25 & 0.02 &0.18 & 30.24 & 0.04 & 0.25& \underline{6.63} & 0.00 &0.20& \underline{6.05} & 0.00 & 0.00 \\ \hline
	\multicolumn{2}{l}{Average}  & 0.033 &0.202  & & 0.088 &0.274  & & 0.060 &0.116 &  & 0.027 & 0.011  \\
	\multicolumn{2}{l}{Maximum}  & 0.09 & 0.38 & & 0.34 & 0.62 &  & 0.27 & 0.20 &  & 0.27 & 0.11  \\ \hline
\end{tabular}
\label{table:experiments1}
\\Notes: Results are based on~$10^6$-period Monte Carlo simulations ($10^4$ burn-in periods). We discretize the exponential demand with a step size of 0.1. The underlined values specify that the optimal policy is CBS.
\end{table}
\begin{figure}[h]
    \centering
   \begin{subfigure}[h]{0.49\textwidth}
     \centering
        \includegraphics[scale=0.5]{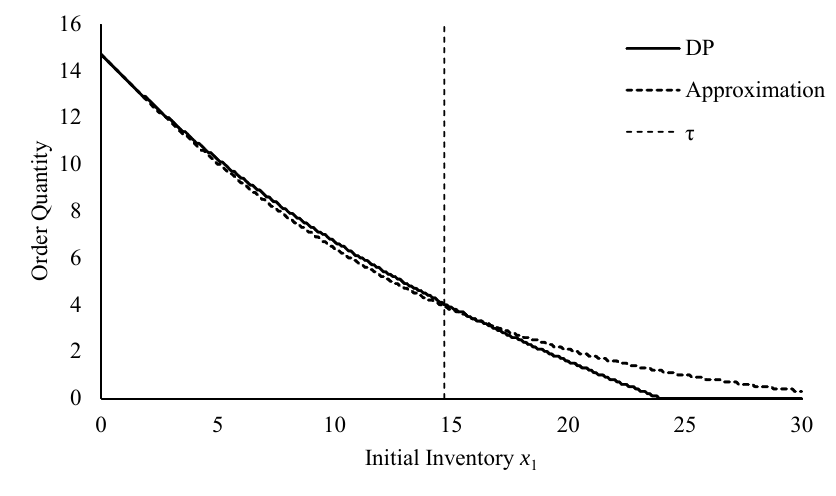}
        \caption{$h=1, m=2$}
        \label{fig:expo10105m2}
    \end{subfigure}
    \begin{subfigure}[h]{0.49\textwidth}
     \centering
        \includegraphics[scale=0.5]{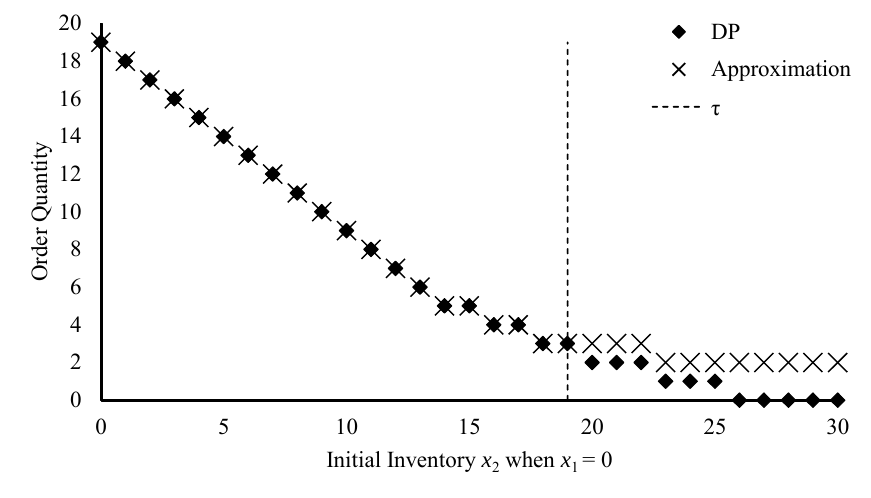}
        \caption{$h=0, m=3$}
        \label{fig:Poi00105m3}
    \end{subfigure}
    \caption{Order amounts; $r=10, \theta =5$. $\tau$ separates the recurrent and non-recurrent inventory levels.}
    \label{fig:ExponentialDemand}
\end{figure}

We also examine the accuracy of the approximate optimal policy with respect to individual order amounts. Figs.~\ref{fig:expo10105m2} and~\ref{fig:Poi00105m3} show examples of this comparison: We observe that the order amounts following policy~$q^{\ast}_h$ closely match those from the optimal policy in recurrent regions;\footnote{The threshold for the recurrent region is specified by the maximum possible initial inventory (i.e.,~$q_h^{\ast}(\mathbf{0})$), as the initial inventory level cannot exceed the order amount at~$\mathbf{x}=0$.} the discrepancies in non-recurrent regions do not impact average costs. In Table~\ref{table:experiments1}, we report the mean absolute deviation (MAD) between the order amounts of the optimal and approximate order policies for the recurrent initial inventory levels.\footnote{For~$m=3$, we compare the order amounts when~$x_1=0$ and~$x_2\geq 0$.}

\section{Concluding Remarks}
We develop a framework to convert a stochastic dynamic inventory model to a single-decision model, by capturing the complex interactions of multi-period decisions in a single externality term. The resulting single-decision model is simpler to solve and yields not only near-optimal average cost, but also close-to-optimal initial inventory-dependent ordering amounts, implying our method captures the fundamental properties of the optimal policy.
Our framework has the potential to be applied to various perishable inventory models, including the more advanced models~(like those considered in~\cite{abouee2019managing} and~\cite{kouki2020analysis}, for example), with appropriate modifications.

\bibliographystyle{model1-num-names}
\bibliography{perishable_Bib}

\newpage
\setcounter{page}{1}
\appendixtitleon
\begin{appendices}{\textbf{\large Online Appendix for ``A Marginal Analysis Framework to Incorporate the Externality
Effect of Ordering Perishables''}}

\section{Deriving Eq.~\eqref{eq:stationary optimality condition}}
We apply variational principle for the functional $L(q)$.
\begin{align*}
\dfrac{\delta L(q)}{\delta q(\mathbf{x})}&= \int_\Omega \left[ \dfrac{\partial L(q(\mathbf{k}),\mathbf{k})}{\partial q(\mathbf{k})} \dfrac{\delta q(\mathbf{k})}{\delta q(\mathbf{x})} f_{\mathbf{X}}^{q}(\mathbf{k}) + L(q(\mathbf{k}),\mathbf{k})\dfrac{\partial f_{\mathbf{X}}^q(\mathbf{k})}{\partial q(\mathbf{x})}\right]d\mathbf{k}\nonumber \\
&= \int_\Omega \left[ \dfrac{\partial L(q(\mathbf{k}),\mathbf{k})}{\partial q(\mathbf{k})} \delta (\mathbf{k}-\mathbf{x}) f_{\mathbf{X}}^{q}(\mathbf{k}) + L(q(\mathbf{k}),\mathbf{k})\dfrac{\partial f_{\mathbf{X}}^q(\mathbf{k})}{\partial q(\mathbf{x})}\right]d\mathbf{k}\nonumber \\
&= \int_\Omega \dfrac{\partial L(q(\mathbf{k}),\mathbf{k})}{\partial q(\mathbf{k})} f_{\mathbf{X}}^{q}(\mathbf{k}) \delta (\mathbf{k}-\mathbf{x}) d\mathbf{k} + \int_\Omega  L(q(\mathbf{k}),\mathbf{k})\dfrac{\partial f_{\mathbf{X}}^q(\mathbf{k})}{\partial q(\mathbf{x})} d\mathbf{k}\nonumber \\
&=\dfrac{\partial L(q(\mathbf{x}),\mathbf{x})}{\partial q(\mathbf{x})} f_{\mathbf{X}}^{q}(\mathbf{x})+ \int_\Omega L(q(\mathbf{k}),\mathbf{k})\dfrac{\partial f_{\mathbf{X}}^q(\mathbf{k})}{\partial q(\mathbf{x})}d\mathbf{k} = 0, \quad \forall \mathbf{x} \in \Omega,
\end{align*}
where we apply \emph{chain rule} first and then \emph{product rule} of the functional derivative \cite[Appendix A of][]{parr1989density} to obtain the first line, replace ${\delta q(\mathbf{k})}/{\delta q(\mathbf{x})}$ with the Dirac delta function $\delta(\mathbf{k}-\mathbf{x})$ to obtain the second line, and apply its sifting property \cite[$g(\mathbf{x})=\int_\Omega g(\mathbf{k}) \delta(\mathbf{k}-\mathbf{x})d\mathbf{k}$ for every continuous function $g(\cdot)$; see, e.g.,][]{bracewell1986fourier} to obtain the fourth line.

\section{Proof of Proposition~\ref{prop:convolution formula}}
We can rewrite~\eqref{eq: Ri} as
\begin{align*}
R_{i}= \max\{{x}^{i}-D^{i}(\mathbf{x}^{i-1}),0\}, \quad i \geq 1,
\end{align*}
which is equivalent to
\begin{align*}
D^{i}(\mathbf{x}^{i-1})+R_{i}= \max\{D^{i}(\mathbf{x}^{i-1}),{x}^{i}\}, \quad i \geq 1.
\end{align*}
Therefore, for~$i \ge 1$,
\begin{align}
F_{D^{i}(\mathbf{x}^{i-1})+R_{i}}(\zeta) =
  \begin{cases}
   F_{D^{i}(\mathbf{x}^{i-1})}(\zeta) & \quad \text{if } \zeta \geq x^i, \nonumber \\
   0  & \quad \text{if } \zeta < x^i.
  \end{cases}
\end{align}
Combining this result with~\eqref{eq: Di}, we obtain
\begin{align*}
F_{D^{i+1} (\mathbf{x}^{i})}(z) &=Pr\{D^{i+1}(\mathbf{x}^{i})\leq z \}=Pr\{D^{i}(\mathbf{x}^{i-1})+R_{i}+D_{i+1} \leq z \}\\
&=Pr\{D^{i}(\mathbf{x}^{i-1})+R_{i} \leq z-D_{i+1} \}= \int_{\xi=-\infty}^{\infty}  F_{D^{i}(\mathbf{x}^{i-1})+R_i} (z-\xi) f_{D_{i+1}}(\xi) d \xi\\
&=\int_{\xi=0}^{z-x^i}  F_{D^{i} (\mathbf{x}^{i-1})} (z-\xi) f_{D_{i+1}}(\xi) d \xi, \quad \text{if } z > x^i, i \geq 1,
\end{align*}
and~$F_{D^{i+1} (\mathbf{x}^{i})}(z) =0$, {} if~$z \le x^i, i \geq 1.$ \hfill~$\square$

\section{Deriving Eq.~\eqref{eq:marginal external cost}}
We discretize the continuous order-up-to level $q_c$ with a step size of $\Delta$.
\begin{align*}
w_{ex}(q_c)&=\int_\Omega \int_{0}^{q_c}(q_c -z)f_{D^m(\mathbf{k})}(z)dz \dfrac{\partial f_{\mathbf{X}}^{q_c}(\mathbf{k})}{\partial q_c} d \mathbf{k}\\
&=\int_\Omega \int_{0}^{q_c}(q_c -z)f_{D^m(\mathbf{k})}(z)dz \dfrac{f_{\mathbf{X}}^{q_c+\Delta}(\mathbf{k})-f_{\mathbf{X}}^{q_c}(\mathbf{k})}{\Delta} d \mathbf{k}\\
&=\dfrac{\int_\Omega \int_{0}^{q_c}(q_c -z)f_{D^m(\mathbf{k})}(z)dz f_{\mathbf{X}}^{q_c+\Delta}(\mathbf{k}) d \mathbf{k}-\int_\Omega \int_{0}^{q_c}(q_c -z)f_{D^m(\mathbf{k})}(z)dz f_{\mathbf{X}}^{q_c}(\mathbf{k}) d \mathbf{k}}{\Delta}\\
&=\frac{\mathbb{E}[q_c - D^m(\mathbf{X}^{q_c+\Delta})]^{+}-\mathbb{E}[q_c -D^m(\mathbf{X}^{q_c})]^{+}}{\Delta}\\
&=\dfrac{n_w^\Delta(q_c)-n_w(q_c)}{\Delta }.
\end{align*}

\section{Proof of Proposition~\ref{prop:external cost}}
We assume (as in \S\ref{sec:application}) that~$h \geq 0$,~$r-c>0$,~$\theta+c>0$, and~$f_D (d)>0, \forall d \ge 0$. Since we discuss the properties of random initial inventory vectors, it is convenient to use the concept of the first-order stochastic dominance (FSD), which is defined as follows:
\begin{definition}
\emph{A random variable~$X$ first-order stochastically dominates another random variable~$Y$ ($X \succeq_{FSD} Y$) if and only if~$F_{X}(x) \leq F_{Y}(x)$,~$\forall x \in \mathbb{R}$.}
\label{def:def2}
\end{definition}
For notational convenience, we write~$\mathbf{X} \succeq_{FSD} \mathbf{Y}$ for random vectors~$\mathbf{X}$ and~$\mathbf{Y}$ if
the FSD property holds componentwise:~$X_i \succeq_{FSD} Y_i$ for all~$i^{th}$ elements of~$\mathbf{X}$ and~$\mathbf{Y}$.
To prove the FSD property for random variables, the following property is convenient and well-known \citep[see, e.g.,][]{Wolfstetter:99}.
\begin{property}
\label{prop:fsd}
$X \succeq_{FSD} Y \iff \mathbb{E}_{X}[f(X)] \geq \mathbb{E}_{Y}[f(Y)]$ for any non-decreasing function~$f(\cdot)$.
\end{property}

Next, we present two useful FSD relationships for~$D^m(\cdot)$. Let~$\mathbf{x}$ and~$\mathbf{y}$ be initial inventory vectors in~$\Omega=\mathbb{R}_{\geq0}^{m-1}$, and~$\mathbf{X}$ and~$\mathbf{Y}$ be the corresponding random vectors.
\begin{lemma}
\label{lem:effective demand1}
$\mathbf{x} \ge \mathbf{y}$ component-wise~$\implies D^m (\mathbf{x}) \succeq_{FSD} D^m (\mathbf{y})$.
\end{lemma}
\noindent \textbf{Proof of Lemma~\ref{lem:effective demand1}:} The proof follows from Proposition~\ref{prop:convolution formula} and Definition~\ref{def:def2}. \hfill~$\square$
\begin{lemma}
\label{lem:effective demand2}
$\mathbf{X} \succeq_{FSD} \mathbf{Y} \implies D^m (\mathbf{X}) \succeq_{FSD} D^m (\mathbf{Y})$.
\end{lemma}
\noindent \textbf{Proof of Lemma~\ref{lem:effective demand2}:} Combining Property~\ref{prop:fsd} and Lemma~\ref{lem:effective demand1}, we have~$\mathbf{x} \ge \mathbf{y}$ componentwise~$\Longrightarrow \mathbb{E}_{D^m (\mathbf{x})}[f(D^m (\mathbf{x}))] \geq \mathbb{E}_{D^m (\mathbf{y})}[f(D^m (\mathbf{y}))]$ for any non-decreasing function~$f(\cdot)$. This result indicates that~$g(\mathbf{x}) \doteq \mathbb{E}_{D^m (\mathbf{x})}[f(D^m (\mathbf{x}))]=\mathbb{E}_{D^m(\mathbf{X})|\mathbf{X}=\mathbf{x}}[f(D^m (\mathbf{X}))|\mathbf{X}=\mathbf{x}]$ is a non-decreasing function in~$\mathbf{x}$ componentwise (because~$g(\mathbf{x}) \ge g(\mathbf{y})$ whenever~$\mathbf{x} \ge \mathbf{y}$ componentwise). Using Property \ref{prop:fsd} once again this time with~$g(\mathbf{x})$ we define above and the law of total expectation,~$\mathbf{X} \succeq_{FSD} \mathbf{Y}\Longrightarrow \mathbb{E}_{\mathbf{X}}[g(\mathbf{X})] \geq \mathbb{E}_{\mathbf{Y}}[g(\mathbf{Y})] \Longleftrightarrow \mathbb{E}_{D^m (\mathbf{X})}[f(D^m (\mathbf{X})] \geq \mathbb{E}_{D^m (\mathbf{Y})}[f(D^m (\mathbf{Y})]$ for any non-decreasing function~$f(\cdot)$, which indicates~$D^m (\mathbf{X}) \succeq_{FSD} D^m (\mathbf{Y})$. \hfill~$\square$

Let~$\mathbf{X}^{q_c}$ and~$X_m^{q_c}$ be the initial inventory random vector and the new order under the CBS policy~$q_c$, respectively. Let~$\tilde{\mathbf{X}}^{q_c}=(\mathbf{X}^{q_c}, X_m^{q_c}) \in \mathbb{R}^m_{\ge0}$. Since the entire inventory follows the CBS policy~$q_c$,~$\sum_{i=1}^{m}X_i^{q_c}=q_c$ must hold. Consider increasing the order-up-to level~$q_c$ by a positive infinitesimal~$\delta_{q_c}$. Then the stationary distribution of the entire inventory (including the new order) shifts from~$\tilde{\mathbf{X}}^{q_c}$ to~$\tilde{\mathbf{X}}^{q_c+\delta_{q_c}}$. The following relationship holds:
\begin{lemma}
\label{lem:initial inventory}
$\mathbf{X}^{q_c+\delta q_c} \succeq_{FSD} \mathbf{X}^{q_c}$.
\end{lemma}
\noindent \textbf{Proof of Lemma~\ref{lem:initial inventory}:} Define a discrete time stochastic process~$\{{\mathbf{\tilde X}}^{q_c}(t), t=0,1,2,...\}$ to represent the entire inventory at time period~$t \in \mathbb{Z}_{\ge0}$. Consider a sample path~${\mathbf{\tilde X}}^{q_c}(t;\omega)$. Without loss of generality, we assume~$\mathbf{X}^{q_c}(0;\omega) = \mathbf{0}$ and~$X_m^{q_c}(0;\omega)=q_c$, which repeatedly appear one period after we encounter a shortage of inventory (note:~$\mathbf{x}=\mathbf{0}$ is recurrent). Suppose the CBS policy is modified from~$q_c$ to~$q_c+\delta_{q_c}$, where~$\delta_{q_c}$ is a positive infinitesimal that is non-divisible. Then the sample path at~$t=0$ shifts from~${\mathbf{\tilde X}}^{q_c}(0;\omega)=(\mathbf{0},q_c)$ to~${\mathbf{\tilde X}}^{q_c+\delta_{q_c}}(0;\omega)=(\mathbf{0},q_c+\delta_{q_c})$. Assuming that this~$\delta_{q_c}$ is used last in each age category, either one of the two occurs every period: (1)~$\delta_{q_c}$ is not used, in which case~$\delta_{q_c}$ becomes older (or wasted) and shows up in the older age category (or the new order category) in the next period, or (2)~$\delta_{q_c}$ is used, in which case~$\delta_{q_c}$ shows up in the same or newer age category in the next period. Hence, the revised sample path is represented as~${\mathbf{\tilde X}}^{q_c+\delta_{q_c}}(t;\omega)={\mathbf{\tilde X}}^{q_c}(t;\omega)+\delta_{q_c}\mathbf{I}(t;\omega),$ where~$\mathbf{I}$ is a random unit vector (one of the age category is 1 and all others are 0) and~$\mathbf{I}(t;\omega)$ is its sample path. It follows that, for each age category~$i$,~$F_{X_i^{q_c+\delta_{q_c}}}(x)=Pr\{X_i^{q_c+\delta_{q_c}} \le x\}=Pr\{X_i^{q_c} +\delta_{q_c}I_i \le x\} \le Pr\{X_i^{q_c} \le x\}=F_{X_i^{q_c}}(x)$,~$\forall x \in [0,\infty)$. Hence, from Definition \ref{def:def2}, we obtain~$\mathbf{\tilde X}^{q_c+\delta q_c} \succeq_{FSD} \mathbf{\tilde X}^{q_c}$, and therefore,~$\mathbf{X}^{q_c+\delta q_c} \succeq_{FSD} \mathbf{X}^{q_c}$. (Note:~$\mathbf{\tilde X}$ and~$\mathbf{I}$ are not independent, but the dependency does not affect the conclusion.) \hfill~$\square$
\begin{lemma}
\label{lem:effective demand3}
$D^m (\mathbf{X}^{q_c+\delta q_c}) \succeq_{FSD} D^m (\mathbf{X}^{q_c})$.
\end{lemma}
\noindent \textbf{Proof of Lemma~\ref{lem:effective demand3}:} The result is immediately obtained from Lemmas~\ref{lem:effective demand2} and~\ref{lem:initial inventory}. \hfill~$\square$

Using Property~\ref{prop:fsd} and Lemma~\ref{lem:effective demand3}, we can bound the externality term and obtain Proposition~\ref{prop:external cost}.

\noindent \textbf{Proof of Proposition~\ref{prop:external cost}:} To prove this property, we rewrite the partial derivative with the expression using a positive infinitesimal change~$\delta_{q_c}$:~$\partial f_{D^m(\mathbf{X}^{q_c})}(z)/\partial q_c=[f_{D^m(\mathbf{X}^{q_c+\delta_{q_c}})}(z) -f_{D^m(\mathbf{X}^{q_c})}(z)] /\delta_{q_c}$.

First part ($w_{ex}(q_c) \le 0$): Changing the order of two integrations and the partial derivative in~\eqref{eq:marginal external cost2}, we obtain:
\begin{align*}
w_{ex}(q_c) &=\int_{0}^{q_c}(q_c -z) \dfrac{\partial \int_\Omega f_{D^m(\mathbf{k})}(z) f_{\mathbf{X}}^{q_c}(\mathbf{k})d\mathbf{k}}{\partial q_c} dz =\int_{0}^{q_c}(q_c -z) \dfrac{\partial f_{D^m(\mathbf{X}^{q_c})}(z)}{\partial q_c} dz \nonumber \\
&=\dfrac{\int_{0}^{q_c}(q_c -z)f_{D^m(\mathbf{X}^{q_c+\delta_{q_c}})}(z)dz - \int_{0}^{q_c}(q_c -z)f_{D^m(\mathbf{X}^{q_c})}(z)dz}{\delta_{q_c}} \nonumber \\
&=- \left( \dfrac{-\mathbb E[q_c-D^m(\mathbf{X}^{q_c+\delta_{q_c}})]^+ + \mathbb E[q_c-D^m(\mathbf{X}^{q_c})]^+}{\delta_{q_c}} \right)\le 0,
\end{align*}
where we apply Property~\ref{prop:fsd} and Lemma~\ref{lem:effective demand3} to~$-[q_c-x]^+$, which is a non-decreasing function of~$x$.

Second part ($w_{ex}(q_c) > -1$): Since~$[q_c - D^m(\mathbf{x})]^+$ represents the amount of wastage,~$\mathbb E[q_c-D^m(\mathbf{X}^{q_c})]^+ = \mathbb E[X_1^{q_c}-D]^+$ should hold. Also, under the assumption~${f_D(d)>0,\forall d \ge0}$, $F_{D^m(\mathbf{X}^{q_c})}(q_c) <1$ for a finite~$q_c$. Using these properties, we obtain:
\begin{align*}
w_{ex}(q_c) &=\int_{0}^{q_c}(q_c -z) \dfrac{\partial f_{D^m(\mathbf{X}^{q_c})}(z)}{\partial q_c} dz =\dfrac{\partial \int_{0}^{q_c}(q_c -z)f_{D^m(\mathbf{X}^{q_c})}(z)dz}{\partial q_c}-\int_{0}^{q_c} f_{D^m(\mathbf{X}^{q_c})}(z) dz\nonumber \\
&=\dfrac{\partial \mathbb E[q_c-D^m(\mathbf{X}^{q_c})]^+}{\partial {q_c}}-F_{D^m(\mathbf{X}^{q_c})}(q_c)=\dfrac{\partial \mathbb E[X_1^{q_c}-D]^+}{\partial {q_c}}-F_{D^m(\mathbf{X}^{q_c})}(q_c) \nonumber \\
&=\dfrac{\mathbb E[X_1^{q_c+\delta_{q_c}}-D]^+ - \mathbb E[X_1^{q_c}-D]^+}{\delta_{q_c}}-F_{D^m(\mathbf{X}^{q_c})}(q_c) >-1,
\end{align*}
where we apply Property~\ref{prop:fsd} and Lemma~\ref{lem:initial inventory} to~$\mathbb E[x-D]^+$, which is a non-decreasing function of~$x$. \hfill~$\square$

\section{Proof of Proposition~\ref{thm:existence}}
We continue to assume (as in \S\ref{sec:application}) that~$h \geq 0$,~$r-c>0$,~$\theta+c>0$, and~$f_D (d)>0, \forall d \ge 0$. Proposition~\ref{thm:existence} follows from Proposition~\ref{prop:external cost}.

Let~$g_{\mathbf{x}}(q) \doteq (\theta+c) F_{D^m(\mathbf{x})}(q)-(h+r-c)\bar{F}_D(q)+h+(\theta+c) w_{ex}(q^{\ast}_c)$, where~$q^{\ast}_c$ is independent of~$q$. Observe that~$g_{\mathbf{x}}(q)$ is an increasing function with respect to~$q$ as~${\partial g_{\mathbf{x}}(q)}/{\partial q}=(\theta+c)f_{D^m(\mathbf{x})}(q)+(h+r-c)f_D(q) > 0$. Also, using Proposition~\ref{prop:external cost},~$g_{\mathbf{x}}(0) =-(r-c)+(\theta+c)w_{ex}(q^{\ast}_c)<0$. Finally there always exists a finite~$\hat{q}$ satisfying~$F_{D^m(\mathbf{x})}(\hat{q})>1-\hat{\delta}$, where~$\hat{\delta} = \frac{(\theta+c)(1+w_{ex}(q^{\ast}_c))}{h+r+\theta} \in (0,1)$. Such~$\hat{q}$ also satisfies~$F_{D}(\hat{q})>1-\hat{\delta}$ (or equivalently,~$\bar{F}_D(\hat{q})<\hat{\delta}$) since~$D^m(\mathbf{x}) \succeq_{FSD} D$ (see~\eqref{eq: Di} and Definition~\ref{def:def2}). Using Proposition~\ref{prop:external cost} and~$\hat{\delta}$ defined above, this~$\hat{q}$ satisfies
\begin{align*}
g_{\mathbf{x}}(\hat{q}) &=(\theta+c) F_{D^m(\mathbf{x})}(\hat{q})-(h+r-c)\bar{F}_D(\hat{q})+h+(\theta+c) w_{ex}(q^{\ast}_c)\\
&>(\theta+c)(1-\hat{\delta})-(h+r-c)\hat{\delta}+h+(\theta+c) w_{ex}(q^{\ast}_c) \\
&=(\theta+c)(1+w_{ex}(q^{\ast}_c))-(h+\theta+r)\hat{\delta}+h =h\ge 0.
\end{align*}
Since~$g_{\mathbf{x}}(q)$ is monotonic, we can conclude that there exists a unique and finite solution~$q_h^\ast(\mathbf{x}) \in (0,1)$ that satisfies~$g_{\mathbf{x}}(q)=0$ (and hence the optimality condition~\eqref{eq:optimality3}) for any initial inventory vector~$\mathbf{x} \in \Omega$. \hfill~$\square$

Finally, we present a corollary of Proposition~\ref{thm:existence}, which is utilized in the proof of Proposition~\ref{thm:CBS condition}.
\begin{corollary}
\label{cor:same q same F}
Consider two initial inventory vectors~$\mathbf{x_1},\mathbf{x_2} \in \Omega$, where~$\mathbf{x_1} \ne \mathbf{x_2}$. Then~$$|q_h^{\ast}(\mathbf{x_1}) - q_h^{\ast}(\mathbf{x_2})| \to 0 \iff |F_{D^m(\mathbf{x_1})}(q_h^{\ast}(\mathbf{x_1})) - F_{D^m(\mathbf{x_2})}(q_h^{\ast}(\mathbf{x_2}))| \to 0.$$
\end{corollary}
\noindent \textbf{Proof of Corollary~\ref{cor:same q same F}:} Since the solution to~\eqref{eq:optimality3} is unique and finite (Proposition~\ref{thm:existence}), we can write two optimality equations corresponding to initial inventory vectors~$\mathbf{x_1}$ and~$\mathbf{x_2}$. Note that~$w_{ex}(q^{\ast}_c)$ in~\eqref{eq:optimality3} does not depend on~$\mathbf{x}$. By subtracting one from the other, we obtain
$$(\theta+c) (F_{D^m(\mathbf{x_1})}(q_h^{\ast}(\mathbf{x_1})) - F_{D^m(\mathbf{x_2})}(q_h^{\ast}(\mathbf{x_2})))=(h+r-c) (\bar{F}_D(q_h^{\ast}(\mathbf{x_1})) - \bar{F}_D(q_h^{\ast}(\mathbf{x_2}))).$$
The result follows from the assumptions~$\theta+c>0$,~$h+r-c>0$, and continuous~$F_D(d)$ for $d \in [0,+\infty)$. \hfill~$\square$

\section{Proof of Proposition~\ref{thm:CBS condition}}
We first show the relationship between two solutions with different initial inventory levels (Lemma~\ref{lem:prop2}), from which we can determine the upper and lower bounds of the solution (Lemma~\ref{lem:bound}). If the gap between the upper and lower bounds shrinks, a state-dependent policy should approach CBS. The condition to make the gap shrink is provided in Proposition~\ref{thm:CBS condition}.
\begin{lemma}
\label{lem:prop2}
$q_h^{\ast}(\mathbf{x_1}) \geq q_h^{\ast}(\mathbf{x_2})$ and~$F_{D^m(\mathbf{x_1})}(q_h^{\ast}(\mathbf{x_1})) \leq F_{D^m(\mathbf{x_2})}(q_h^{\ast}(\mathbf{x_2}))$ if~$\mathbf{x_1} \ge \mathbf{x_2}$ component-wise.
\end{lemma}

\noindent \textbf{Proof of Lemma~\ref{lem:prop2}:} As in the proof of Proposition~\ref{thm:existence}, we define~$g_{\mathbf{x}}(q) \doteq (\theta+c) F_{D^m(\mathbf{x})}(q)-(h+r-c)\bar{F}_D(q)+h+(\theta+c) w_{ex}(q^{\ast}_c)$.  This~$g_{\mathbf{x}}(q)$ is an increasing function with respect to~$q$. Now, let~$q_h^{\ast}(\mathbf{x_1})$ and~$q_h^{\ast}(\mathbf{x_2})$ be the unique, finite solutions to~$g_{\mathbf{x_1}}(q)=0$ and~$g_{\mathbf{x_2}}(q)=0$, respectively. Since~$\mathbf{x_1} \ge \mathbf{x_2}$ componentwise~$\Longrightarrow D^m (\mathbf{x_1}) \succeq_{FSD} D^m (\mathbf{x_2})$ (Lemma~\ref{lem:effective demand1})~$\iff F_{D^m (\mathbf{x_1})}(q) \leq F_{D^m (\mathbf{x_2})}(q)$,~$\forall q \in \mathbb{R}$ (Definition \ref{def:def2}), it follows that~$g_{\mathbf{x_1}}(q) \leq g_{\mathbf{x_2}}(q), \forall q \in \mathbb{R}$. In particular, at~$q=q_h^{\ast}(\mathbf{x_2})$, we obtain~$g_{\mathbf{x_1}}(q_h^{\ast}(\mathbf{x_2})) \leq g_{\mathbf{x_2}}(q_h^{\ast}(\mathbf{x_2})) =0$, which implies~$q_h^{\ast}(\mathbf{x_1}) \geq q_h^{\ast}(\mathbf{x_2})$. Furthermore,~$q_h^{\ast}(\mathbf{x_1}) \geq q_h^{\ast}(\mathbf{x_2})$ implies~$\bar{F}_D(q_h^{\ast}(\mathbf{x_1})) \leq \bar{F}_D(q_h^{\ast}(\mathbf{x_2}))$ because~$\bar{F}_D(q)$ is a decreasing function of~$q$. Combining this result with~\eqref{eq:optimality3}, we obtain~$F_{D^m(\mathbf{x_1})}(q_h^{\ast}(\mathbf{x_1})) \leq F_{D^m(\mathbf{x_2})}(q_h^{\ast}(\mathbf{x_2}))$. \hfill~$\square$

\begin{lemma}
\label{lem:bound}
$q^{\dagger} \leq q_h^{\ast}(\mathbf{x}) \leq q^{\ddagger}$ and $F_{D^m}(q^{\dagger}) \geq F_{D^m(\mathbf{x})}(q_h^{\ast}(\mathbf{x}))$, $\forall \mathbf{x} \in \Omega$.
\end{lemma}

\noindent \textbf{Proof of Lemma~\ref{lem:bound}:} From Lemma~\ref{lem:prop2}, we have~$q_h^{\ast}(\mathbf{0}) \leq q_h^{\ast}(\mathbf{x}) \leq q_h^{\ast}(\mathbf{y})$ and~$F_{D^m(\mathbf{0})}(q_h^{\ast}(\mathbf{0})) \geq F_{D^m(\mathbf{x})}(q_h^{\ast}(\mathbf{x})), \forall \mathbf{y} \ge \mathbf{x} (\in \Omega)$ componentwise. We obtain the result by taking the limit of a large initial inventory~$\mathbf{y}$ and denoting~$D^m=D^m(\mathbf{0})$,~$q^\dagger=q_h^{\ast}(\mathbf{0})$, and~$q^{\ddagger}=\lim_{v \to \infty} q_h^{\ast}(\mathbf{y})$, where~$v$ is the smallest component in the initial inventory vector~$\mathbf{y}$.  \hfill~$\square$

\noindent \textbf{Proof of Proposition~\ref{thm:CBS condition}:} We split the proof in three parts:

\noindent (First part:~$F_{D^m}(q^{\dagger}) \to 0 \implies |q^{\ddagger} - q^{\dagger}| \to 0$) Using Lemma~\ref{lem:bound},~$F_{D^m}(q^{\dagger}) \to 0 \implies F_{D^m(\mathbf{x})}(q_h^{\ast}(\mathbf{x})) \to 0, \forall \mathbf{x} \in \Omega$. Since~$F_{D^m(\mathbf{x})}(q_h^{\ast}(\mathbf{x}))$ converges to the same value (0) for any initial inventory vector~$\mathbf{x}$, using Corollary~\ref{cor:same q same F}, we can conclude~$|q^{\ddagger} - q^{\dagger}| \to 0$.

\noindent (Second part:~$|q^{\ddagger} - q^{\dagger}| \to 0 \implies q_h^{\ast}(\mathbf{x}) \to q_c^{\ast}, \forall \mathbf{x} \in \Omega_r$) This part is trivial because~$\Omega_r \subseteq \Omega$.

\noindent (Third part:~$q_h^{\ast}(\mathbf{x}) \to q_c^{\ast}, \forall \mathbf{x} \in \Omega_r \implies F_{D^m}(q^{\dagger}) \to 0$) Consider two initial inventory vectors:~$\mathbf{x_1}=\mathbf{0}$ and~$\mathbf{x_2}=(q^\dagger,0,\dots,0)$;~$\mathbf{x_2}$ represents~$q^\dagger (= q_h^{\ast}(\mathbf{0}) = q_h^{\ast}(\mathbf{x_1}))$ units of initial inventory with remaining lifetime of~$m-1$ periods. Note that~$\mathbf{x_1} \ne \mathbf{x_2}$ and both~$\mathbf{x_1},\mathbf{x_2} \in \Omega_r$ because we assume that~$D$ can take 0 and any large amount. Note also that from Proposition~\ref{prop:convolution formula}, we know~$F_{D^m(\mathbf{x_2})}(q^\dagger)=0$. (This is intuitively obvious:~$D^m(\mathbf{x_2})$ is the total outflow (through demand and wastage) from periods 1 to~$m$ (excluding the wastage in period~$m$) when the initial inventory is~$\mathbf{x_2}$. Hence, the support of its CDF is bounded below by~$q^\dagger$.) Now, suppose~$q_h^{\ast}(\mathbf{x}) \to q_c^{\ast}, \forall \mathbf{x} \in \Omega_r$, then~$q_h^{\ast}(\mathbf{x_2}) \to q_h^{\ast}(\mathbf{x_1}) = q^\dagger$. Therefore, using Corollary~\ref{cor:same q same F} and replacing~$q_h^{\ast}(\mathbf{x_2})$ with~$q^\dagger$, we obtain~$F_{D^m}(q^{\dagger})=F_{D^m(\mathbf{x_1})}(q_h^{\ast}(\mathbf{x_1})) \to F_{D^m(\mathbf{x_2})}(q_h^{\ast}(\mathbf{x_2})) \to F_{D^m(\mathbf{x_2})}(q^\dagger)=0.$ \hfill~$\square$

\end{appendices}

\end{document}